\documentclass[conference]{IEEEtran}
\IEEEoverridecommandlockouts
\usepackage{cite}
\usepackage{amsmath,amssymb,amsfonts}
\usepackage{algorithmic}
\usepackage{graphicx}
\usepackage{textcomp}
\usepackage{xcolor}
\usepackage{array}
\usepackage{booktabs}
\usepackage{hyperref}
\usepackage{multirow}
\usepackage{xcolor}

\DeclareMathOperator*{\argmin}{arg\,min}

\def\BibTeX{{\rm B\kern-.05em{\sc i\kern-.025em b}\kern-.08em
    T\kern-.1667em\lower.7ex\hbox{E}\kern-.125emX}}

\usepackage{enumitem} 

\renewcommand*\descriptionlabel[1]{\hspace\labelsep \normalfont #1}
\makeatother

\begin{document}

\title{Contracting Strategies for Electrolyzers to Secure Grid Connection: The Dutch Case
}

\author{\IEEEauthorblockN{Thomas Swarts\textsuperscript{a}, Jalal Kazempour\textsuperscript{b}, Wouter van den Akker\textsuperscript{a}, Johan Morren\textsuperscript{a}, Arjan van Voorden\textsuperscript{c}, Han Slootweg\textsuperscript{a}}
\IEEEauthorblockA{\textsuperscript{a} Eindhoven University of Technology, \textsuperscript{c} Delft University of Technology, The Netherlands \\
\textsuperscript{b }Department of Wind and Energy Systems, Technical University of Denmark, Kgs. Lyngby, Denmark \\
\textsuperscript{a} \{l.j.m.m.t.swarts, j.morren, w.f.v.d.akker, j.g.slootweg\}@tue.nl, \textsuperscript{b} jalal@dtu.dk, \textsuperscript{c} a.m.vanvoorden-1@tudelft.nl}
}

\maketitle

\begin{abstract}
In response to increasing grid congestion in the Netherlands, non-firm connection and transport agreements (CTAs) and capacity restriction contracts (CRCs) have been introduced, allowing consumer curtailment in exchange for grid tariff discounts or per-MW compensations. This study examines the interaction between an electrolyzer project, facing sizing and contracting decisions, and a network operator, responsible for contract activations and determining grid connection capacity, under the new Dutch regulations. The interaction is modeled using two bilevel optimization problems with alternating leader-follower roles. Results highlight a trade-off between CRC income and non-firm CTA tariff discounts, showing that voluntary congestion management by the network operator increases electrolyzer profitability at CRC prices below \texteuro10/MW but reduces it at higher prices. Furthermore, the network operator benefits more from reacting to the electrolyzer owner’s CTA decisions than from leading the interaction at CRC prices above \texteuro10/MW. Ignoring the other party’s optimization problem overestimates profits for both the network operator and the electrolyzer owner, emphasizing the importance of coordinated decision-making.
%
\end{abstract}

\vspace{2mm}

\vspace{2mm}
\section*{List of Abbreviations}

\begin{description}[labelwidth=1cm, labelsep=1cm]
    \item[CTA] Connection and transport agreement
    \item[CRC] Capacity restriction contract
    \item[FA] Firm CTA  
    \item[NFA] Non-firm CTA  
    \item[NFA85] Non-firm CTA-85
\end{description}

\vspace{1mm}
\section{Introduction}
The growing penetration of distributed energy resources, coupled with the rapid electrification of transport, residential, and industrial sectors, has led to severe congestion challenges in power distribution and transmission grids across Europe. Addressing these congestion issues through grid reinforcements is both time-consuming and expensive, resulting in extensive delays and long waiting lists for new connection and transport agreements (CTAs). The situation is particularly critical in the Netherlands, where congestion affects most medium- and high-voltage grids, with waiting lists in some regions surpassing a hundred pending requests \cite{Netbeheer2025}. 

The inability of consumers to connect to Dutch power grids significantly hinders the progress of the energy transition. Grid congestion has already stalled several green hydrogen projects in the Netherlands \cite{Zoelen2024}. These delays are particularly concerning given the country's ambitious carbon neutrality goals, which require the development of 16 to 45 GW of domestic electrolysis capacity by 2050 \cite{Energiesysteem2023}.

To address these challenges, the Dutch regulator introduced new contract types to the grid code between 2022 and 2024, allowing grid operators to temporarily limit the contracted capacity of consumers to mitigate anticipated congestion.
In 2022, the Dutch regulator implemented capacity restriction contracts (CRCs), followed by non-firm CTAs in 2024 \cite{ACMCode2025}. Under non-firm CTAs, consumers receive flexible (non-firm) transport capacity instead of guaranteed (firm) capacity, in exchange for discounts on their annual transport tariffs. Eligibility for CRCs depends on the type of CTA chosen by the consumer. Consumers with CRCs are compensated by the network operator for each MW of reduced power from their contracted capacity.

Existing studies demonstrate the potential of non-firm CTAs and congestion management in addressing grid congestion and facilitating renewable energy development. For example, \cite{Verhoeven2024} explored budget allocation strategies for redispatch and CRCs, while \cite{Shen2018} proposed a sequential congestion management scheme involving dynamic tariffs and network reconfiguration. Studies on non-firm CTAs include their application to small-scale photovoltaic projects in France \cite{Muller2023}, and their role in maximizing grid access in the Netherlands \cite{Mehmood2024}.
Despite these insights, the literature tends to treat non-firm CTAs and CRCs in isolation, overlooking their interplay and the interactions they create between consumers and network operators. This gap is particularly relevant in the Dutch regulatory context, where consumers can choose from multiple non-firm CTAs and non-firm CTA-CRC combinations.

To address this gap, this study analyzes optimal non-firm CTA and CRC contracting decisions between a client and network operator. The study considers the case of an electrolyzer project in the Netherlands that cannot expand due to grid congestion, a scenario mirrored in real-world projects such as GROHW and Hessenpoort \cite{Zoelen2024}. The electrolyzer owner selects CTAs, while the network operator decides on connecting capacity and activation timing for non-firm CTAs and CRCs. The interaction is modeled using two bilevel optimization problems with alternating leader-follower roles.


The remainder of the paper is organized as follows: Section \ref{sec: contract info} elaborates on the new CTA and CRC contracts in the Dutch grid code. Section \ref{sec: model formulation} develops the mathematical framework. Section \ref{sec: results} presents numerical results, and finally Section \ref{sec: con} concludes the paper.

 \vspace{1mm}
\section{Novel Contracts in the Dutch Power System} \label{sec: contract info}
European network operators, as natural monopolies, are heavily regulated to ensure fair and non-discriminatory consumer access. Traditionally, network operators were required to (\textit{i}) unconditionally provide consumers with their requested capacity on an asset within technical limits, and (\textit{ii}) guarantee access to that capacity at all times. However, with the introduction of the new CTAs and CRCs, these obligations no longer apply in congested areas of the Netherlands\footnote{An area is considered congested when the required transport capacity exceeds the existing grid capacity, which may result from new connection requests, expansions by existing consumers, or changes in market conditions.}.

\vspace{1mm}
\subsection{Non-firm Power Connection and Transport Agreements}
Starting in April 2025, network operators in the Netherlands will offer consumers the choice between three types of CTAs. The firm CTA (FA) guarantees full access to the contracted connection capacity at all times. The non-firm CTA-85 (NFA85) allows the network operator to reduce the consumer’s connection capacity within an annual time budget $B^{\rm{NFA85}}$, limiting capacity for up to 15\% of the hours in a year. Lastly, the non-firm CTA (NFA) enables the network operator to limit connection capacity at any time. However, due to a lack of public acceptance, network operators, in consultation with the regulator, are developing a new NFA that guarantees the consumer a predetermined daily energy budget, $B^{\rm{NFA}}$. This new NFA is assumed in this study.

Consumers who agree to non-firm contracts benefit from discounts on their fixed connection tariffs. They can choose a single contract or split their connection capacity across a combination of contracts (e.g., 1-MW FA, 1-MW NFA85, and 1-MW NFA). To assist consumers in making informed decisions between CTAs, the network operator provides current and simulated future network profiles, along with anticipated congestion. This allows consumers to assess how much transport capacity they can access with non-firm contracts and when. The network operator informs consumers about capacity restrictions by the day before, no later than the closing of the day-ahead (DA) market. Details of the CTA contracts are summarized in  Table \ref{tab: contracs}.

\begin{table}[t]
\centering
\caption{Dutch grid code for connection and transport agreements}
\scriptsize 
\begin{tabular}{l|ccc}
\toprule
& \textbf{FA} (firm) & \textbf{NFA85} (non-firm) & \textbf{NFA} (non-firm) \\
\midrule
Remuneration & None & 100\% discount kW-c & To be determined \\
Transport rights & All times & 85\% of year & 0\% of year  \\
Eligible for CRC & Yes & Yes & No \\
Communication & N/A & Before closure DA & Before closure DA \\
Budget & N/A & Time budget & Energy budget \\
\bottomrule
\end{tabular}
\label{tab: contracs}
\end{table}
\subsection{Capacity Restriction Contracts}

Since 2022, the Dutch grid code includes two instruments for congestion management: redispatch and CRCs. This study focuses exclusively on CRCs. According to the grid code, before network operators can deny new connection requests in a congested network, they must allocate a minimum budget \( B^{\rm{CM}} \) of €1.02/MW of transportable capacity for congestion management, applicable to both existing and new clients in the network. Network operators may also voluntarily allocate a budget beyond this minimum for further congestion management efforts.

Consumers with a capacity of over 1 MW are required to participate in congestion management if they are located in congested areas and their CTA permits it. Similar to non-firm connection and transport agreements, CRCs enable network operators to limit a consumer's connection capacity during peak periods. Consumers are compensated per MW of reduced power from their contracted capacity. The CRC price \( \lambda^{\rm{CRC}} \) can be either fixed or dynamic and is negotiated between the consumer and the network operator. This study assumes that connection capacity can be fully curtailed under CRCs, though specific limits may be defined. Network operators activate CRCs on a day-ahead basis, prior to the closure of the day-ahead market.

Consumers with FA or NFA85 contracts are eligible for CRCs, while consumers with NFA contracts are exempt. Since network operators are obligated to spend \( B^{\rm{CM}} \) on congestion management, this study assumes that FA and NFA85 contracts automatically include a CRC. Additionally, network operators may establish an optional CRC+ contract for congestion management beyond the minimum budget. For consumers with NFA85 contracts, network operators can implement curtailments by activating either CRCs or the NFA85 mechanism, considering the respective budgets \( B^{\rm{NFA85}} \) and \( B^{\rm{CM}} \). Curtailments under NFA85 and CRCs can occur simultaneously, but the same capacity cannot be curtailed twice.


\vspace{2mm}

\section{Model Formulation} \label{sec: model formulation}
In the following, Greek letters and uppercase symbols represent parameters, while lowercase symbols denote variables. We define the following indices and sets: The index for hours is denoted by $t \in \mathcal{T}$. Additionally, $c \in \mathcal{C}$ represents the index for CTAs, where $\mathcal{C} = \{\mathrm{FA}, \mathrm{NFA85}, \mathrm{NFA}\}$. Symbols indexed by $c$, such as $r_{c,t}$, will be used later, where $c=1$, $c=2$ and $c=3$ correspond to the FA, NFA85, and NFA contracts, respectively.

\begin{figure}[t]
    \centering
\includegraphics[trim={0mm 5mm 0mm 4mm},clip,width = 0.8\linewidth]{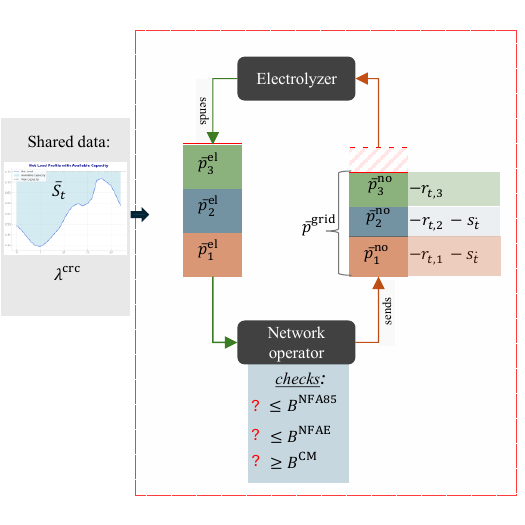}
    \caption{Interaction between the electrolyzer owner and the network operator, along with their respective decisions. In the figure, $s^{.}_{t}$ denotes the CRC curtailments $s_{t}$ and $s^{\rm{+}}_{t}$.}
    \label{fig:interactions}
    \vspace{-2mm}
\end{figure} 

\vspace{2mm}

\subsection{Interaction Between the Electrolyzer Owner and the Network Operator}

We examine a scenario inspired by the GROWH project, where the project owners' request for a connection expansion led to congestion and could not be fulfilled. The project aims to expand to meet a maximum hydrogen demand \( D_{t} \) in the region. The introduction of new CTA and CRC contract types enables the electrolyzer owner to request a larger connection. 

When making contracting decisions, the network operator informs the electrolyzer owner about the residual network capacity and available transport capacity for each hour $t$ within a given time horizon. The electrolyzer owner then submits a connection expansion request to the network operator, specifying the desired capacity for their preferred non-firm CTAs, $\bar{p}^{\rm{el}}_{c}$. The network operator evaluates the request by checking the residual grid capacity, $\bar{S_{t}}$, ensuring compliance with the available CTA budgets ($B^{\rm{NFA85}}$, $B^{\rm{NFA}}$), and verifying the congestion management budget, $B^{\rm{CM}}$. If $B^{\rm{CM}}$ is exceeded, the network operator assesses whether additional congestion costs can be accepted. If the requested capacity exceeds the available limits, the network operator may adjust the capacity to stay within acceptable bounds. Finally, the network operator communicates the available connection capacity, $\bar{p}^{\rm{grid}}$, the contract capacities that can be accommodated, $\bar{p}^{\rm{no}}_{c}$, and the expected curtailments for both CTA contracts ($r_{2,t}$ and $r_{3,t}$) and CRC contracts ($s_{t}$ and $s^{\rm{+}}_{t}$) to the electrolyzer owner. The interaction between the electrolyzer owner and the network operator is illustrated in Figure \ref{fig:interactions}.

The interaction between the electrolyzer owner and the network operator creates a Stackelberg follower-leader game, where the leader makes an initial decision that influences the follower's actions, which in turn affects the leader's outcome. We consider two possible leader-follower hierarchies: 
\begin{enumerate}
    \item [(\textit{i})] In the first hierarchy, the electrolyzer owner leads the game by selecting preferred CTA capacities, anticipating the network operator's decisions on maximum connection capacity and curtailment under CTA and CRC contracts, assuming perfect knowledge of the network operator's optimization problem.
    \item [(\textit{ii})] In the second hierarchy, the network operator leads the game by providing the residual network space and grid connection capacity, anticipating the electrolyzer owner's decisions on contracting capacities, with the assumption of perfect knowledge of the electrolyzer owner's optimization problem.
\end{enumerate}

This study formulates both games as bilevel programs. A bilevel program consists of an upper-level optimization problem (representing the leader's problem) constrained by a lower-level optimization problem (representing the follower's problem). After appropriate reformulations, both bilevel programs ultimately result in mixed-integer linear programs.



\vspace{1mm}
\subsection{Game I: Electrolyzer Owner as the Leader}
We begin with the game in which the electrolyzer owner takes the leading role, while the network operator is the follower. 
The upper-level problem of the electrolyzer owner is formulated in \eqref{1ely_obj_v1}-\eqref{1ely_upper_h2_demand}, subject to the lower-level problem of the network operator, given in  \eqref{1DSO_obj}-\eqref{total curtailment constraint}. The set of upper-level variables, determined by the electrolyzer owner, is $ \Theta^{\rm{el,up}}$ $=\{\bar{p}^{\rm{el}}_{c}$, $ p^{e}_{t}$, $f_{t}\} $, while the set of lower-level variables, determined by the network operator, is $ \Theta^{\rm{no,low}} =$ $\{\bar{p}^{\rm{grid}}$, $ \bar{p}^{\rm{no}}_{c}$, $r_{3, t}$, $r_{2, t}$, $s^{+}_{t}$, $s_{t} \} $. These variables will be defined throughout the text as needed.

The upper-level objective function minimizes the annualized costs minus revenues of the electrolyzer:
%
%
\begin{subequations}
\begin{align}  \label{1ely_obj_v1}
    \min_{\Theta^{\rm{el,up}}} \ \; &C^{\rm{el}}\bar{p}^{\rm{grid}} + \sum_{c\in{\mathcal{C}}}T_{c}\bar{p}^{\rm{el}}_{c} + \sum_{t\in{\mathcal{T}}}\left(\lambda_{t}^{\rm{e}}p_{t}^{\rm{e}} - \lambda^{\rm{H_{2}}}f_{t} \right) \nonumber \\ &- \sum_{t\in{\mathcal{T}}}\left( \lambda^{\rm{crc}}s_{t} + \lambda^{\rm{crc+}}s^{\rm{+}}_{t}  \right), 
\end{align}
where the parameter $C^{\rm{el}}$ represents the annualized capital cost of the electrolyzer plant, and the variable $\bar{p}^{\rm{grid}}$ denotes the grid connection capacity, which is equivalent to the electrolyzer's capacity. In the second term, $T_{c}$ refers to the annual connection and transport tariffs for CTA $c$, while $\bar{p}^{\rm{el}}_{c}$ represents the contracted CTA capacity. In the third term, the electrolyzer consumes power $p^{\rm{e}}_{t}$ at the electricity price $\lambda^{\rm{e}}_{t}$ during hour $t$. The electrolyzer owner sells hydrogen, $f_{t}$, at a fixed hydrogen price, $\lambda^{\rm{H_{2}}}$. In the final term, the electrolyzer owner receives compensation ($\lambda^{\rm{crc}}$ and $\lambda^{\rm{crc+}}$) per MW of curtailed power ($s_{t}$ and $s^{\rm{+}}_{t}$) under CRC and CRC+ contracts  in hour \( t \).

The total contracted CTA capacity is limited by the available grid connection capacity as follows: 
\begin{equation}
\sum_{c\in{\mathcal{C}}}\bar{p}^{\rm{el}}_{c} \leq \bar{p}^{\rm{grid}}. \label{ely_total_contract_constraints}
\end{equation}

The power consumption of the electrolyzer $p_{t}^{\rm{e}}$ in hour $t$ is subject to an upper bound, given by:
\begin{align}
    & p_{t}^{\rm{e}} \leq  - s_{t} - s^{\mathrm{+}}_{t} - r_{2, t} - r_{3, t} + \sum_{c\in{\mathcal{C}}}\bar{p}^{\rm{el}}_{c}, \quad \forall t \in \mathcal{T},\label{ely_power_constraint} 
\end{align}
where $r_{2, t}$ and $r_{3, t}$ are curtailments under the NFA85 and NFA contracts, respectively. In addition, the consumption $p_{t}^{\rm{e}}$ is lower-bounded by: 
\begin{align}
    & p_{t}^{\rm{e}} \geq \alpha^{\rm{min}}\bar{p}^{\rm{grid}}, \quad \forall t \in \mathcal{T},\label{ely_min_load_constraint}
\end{align}
where $\alpha^{\rm{min}}$ is the electrolyzer's minimum load ratio.
The hydrogen production $f_{t}$ of the electrolyzer is constrained by:
\begin{align}
     & f_{t} \leq \eta^{\rm{sys}} p_{t}^{\rm{e}}, \quad \forall t \in \mathcal{T}, \label{ely_production_constraint} \\
      & f_{t} \leq D_{t} \quad \forall t \in \mathcal{T}, \label{1ely_upper_h2_demand}
\end{align}
where \( \eta^{\rm{sys}} \) is the system efficiency (kg/MWh) of the electrolyzer plant and \( D_{t} \) is the maximum hydrogen offtake in the plant's region in hour \( t \).
\end{subequations}

Now, we present the lower-level problem of the network operator. The lower-level objective function minimizes the annual cost of the operator as:
\begin{subequations}
\begin{align} \label{1DSO_obj}
    \min_{\Theta^{\rm{no,low}}} \ \; &\sum_{t\in\mathcal{T}}\left(\lambda^{\rm{crc}}_{t}s_{t} +\lambda^{\rm{crc+}}_{t}s^{\rm{+}}_{t}\right)  - \sum_{c\in{\mathcal{C}}}T_{c}\bar{p}^{\rm{no}}_{c} \nonumber \\ 
    & + \pi\sum_{t \in \mathcal{T}} \left(r_{2, t} + r_{3, t} \right),
\end{align}
where the first term represents the CRC contracting cost, compensating the electrolyzer owner. The second term includes the tariff income of the network operator, with the variable \( \bar{p}^{\rm{no}}_{c} \) denoting the contracted capacity for CTA \( c \). The third term is a small penalty, weighted by the parameter \( \pi \), designed to steer the solution away from CTA curtailment that exceeds the amount necessary to prevent congestion. 

The variable $\bar{p}^{\rm{no}}_{c}$ is constrained by: 
\begin{equation} \bar{p}^{\rm{no}}_{c} \leq \bar{p}^{\rm{el}}_{c}, \quad \forall c \in \mathcal{C},\label{DSO_contract_constraint} \end{equation}
which ensures that the electrolyzer owner in the upper level, and not the network operator, decides on the CTA capacities.

The network operator determines the grid connection capacity of the electrolyzer as: \begin{equation} \sum_{c\in{\mathcal{C}}}\bar{p}^{\rm{no}}_{c} = \bar{p}^{\rm{grid}}. \label{DSO total contract capacity} \end{equation}

The network operator determines the transport capacity in hour $t$ to stay within the residual capacity of the network as:
\begin{equation}
    \bar{p}^{\rm{grid}} - s_{t} - s^{\mathrm{+}}_{t} - r_{2, t} - r_{3, t}  \leq \bar{S}_{t}, \quad \forall t\in \mathcal{T}, \label{grid capacity}
\end{equation}
where \( \bar{S}_{t}\) is the residual grid capacity in hour $t$.

The network operator's time-budgeted curtailment decisions under the NFA85 contract are described by: 
\begin{align}
     & r_{2, t}  \leq \bar{p}^{\rm{no}}_{2}, \quad \forall t \in \mathcal{T}, \\ 
        & r_{2, t} \leq b_{t} \frac{D_{t}}{\eta^{\rm{sys}}}, \quad \forall t \in \mathcal{T}, \\ 
        & \sum_{t \in \mathcal{T}} b_{t} \leq B^{\rm{NFA85}}  |\mathcal{T}|, \\ 
        & b_{t} \in \{0,1\}, \quad \forall t \in \mathcal{T}, \label{eq: binary in lower level}
\end{align}
where the variable $\bar{p}^{\rm{no}}_{2}$ represents the NFA85 contract capacity accommodated by the network operator, $b_{t}$ is a binary variable indicating curtailment activations for the NFA85 contract (with a value of one indicating the activation of curtailment and zero otherwise), $|\mathcal{T}|$ denotes the cardinality of the set $\mathcal{T}$, i.e., the number of hours in the given time horizon, and $B^{\rm{NFA85}} |\mathcal{T}|$ is the CTA compliance budget, which enforces an upper bound on the number of curtailment events over $|\mathcal{T}|$. 

Similarly, the network operator's energy-budgeted curtailment decisions under the NFA contract are described by: 
\begin{align}
   &   r_{3, t} \leq \bar{p}^{\rm{no}}_{3}, \quad \forall t \in \mathcal{T},  \\
   & \sum_{t \in \mathcal{T}} r_{3, t}  \leq B^{\rm{NFA}} \bar{p}^{\rm{no}}_{3}, \label{DSO NFA budget}
\end{align}
where the variable \(\bar{p}^{\rm{no}}_{3}\) is the NFA contract capacity accommodated by the network operator, and \( B^{\rm{NFA}} \) is the corresponding budget.

The network operator must allocate a budget for congestion management to support the operation of the electrolyzer:
\begin{equation}
\lambda^{\rm{crc+}}s^{\rm{+}}_{t} = \theta B^{\rm{CM}}, \quad \forall t \in \mathcal{T},
\end{equation}
where \( B^{\rm{CM}} \) is the budget that the network operator is obliged to allocate for congestion management in the area, distributed among all flexible consumers in the congested network. The parameter \( \theta \) represents the proportion of the budget not allocated to other flexible consumers in the area and is thus exclusively spent on the electrolyzer.

Recall that only FA and NFA85 contracts can be combined with CRC contracts, which is enforced by: 
\begin{equation}
   s_{t} + s^{\mathrm{+}}_{t} + r_{2, t} \leq \bar{p}^{\rm{no}}_{1} + \bar{p}^{\rm{no}}_{2}, \quad \forall t \in \mathcal{T}, \label{total curtailment constraint}
\end{equation}
ensuring that the same capacity is not curtailed twice.
\end{subequations}

The resulting bilevel program for this game, with the electrolyzer owner as the leader and the network operator as the follower, is formulated as:
\begin{subequations} \label{bilevel EN}
\begin{align}
    & \min_{\Theta^{\rm{el,up}}} \quad \eqref{1ely_obj_v1} \label{bilevel_EN1} \\ 
    & \rm{s.t.} \quad \eqref{ely_total_contract_constraints}-\eqref{1ely_upper_h2_demand} \label{bilevel_EN2} \\
    & \bar{p}^{\rm{grid}}, r_{3, t}, r_{2, t}, s^{+}_{t}, s_{t} \in \nonumber \\
    & \quad\quad\quad\argmin_{\Theta^{\rm{no,low}}} \{\eqref{1DSO_obj} \;\: \rm{s.t.} \;\: \eqref{DSO_contract_constraint}-\eqref{total curtailment constraint}\}. \label{bilevel_EN3}
\end{align}
\end{subequations}

\vspace{2mm}
\subsection{Game II: Network Operator as the Leader}
\begin{subequations} 
The upper-level problem for the network operator as the leader in the second game  is the same as in  \eqref{1DSO_obj}-\eqref{total curtailment constraint} in the first game. However, it is subject to the lower-level problem of the electrolyzer (follower), given by \eqref{2ely_obj_v1} constrained by \eqref{ely_total_contract_constraints}, \eqref{ely_min_load_constraint}-\eqref{1ely_upper_h2_demand}, and \eqref{ccc}-\eqref{eq: ely upper nfa85}. The set of upper-level variables, determined by the network operator, is \( \Theta^{\rm{no,up}} = \{\bar{p}^{\rm{grid}}, \bar{p}^{\rm{no}}_{c}, r_{3, t}, r_{2, t}, s^{+}_{t}, s_{t} \} \). Furthermore, the set of lower-level variables, determined by the electrolyzer owner, is  $\Theta^{\rm{el,low}} = \{\bar{p}^{\rm{el}}_c, p^{\rm{e}}_{t}, f_{t}, s^{\rm{el}}_{t}, s^{\rm{el+}}_{t}, r^{\rm{el}}_{2,t}, r^{\rm{el}}_{3,t}\}$.

Formulating the lower-level problem of the electrolyzer owner in the second game requires certain adaptations of the upper-level problem in the first game. The electrolyzer owner's objective  is to minimize the annualized cost (i.e., minus profit) as:
\begin{align}  \label{2ely_obj_v1}
    \min_{\Theta^{\rm{el,low}}} \; &C^{\rm{el}}\bar{p}^{\rm{grid}} + \sum_{c\in{\mathcal{C}}}T_{c}\bar{p}^{\rm{el}}_{c} + \sum_{t\in{\mathcal{T}}}\left(\lambda_{t}^{\rm{e}}p_{t}^{\rm{e}} - \lambda^{\rm{H_{2}}}f_{t} \right) \nonumber \\ &- \sum_{t\in{\mathcal{T}}}\left( \lambda^{\rm{crc}}s^{\rm{el}}_{t} + \lambda^{\rm{crc+}}s^{\rm{el+}}_{t}  \right),
\end{align}
where the additional lower-level CRC curtailment decisions \( s^{\rm{el}}_{t} \) and \( s^{\rm{el+}}_{t} \) enable the lower-level problem to account for revenue through CRC activations. 

The following constraints restrict the electrolyzer owner's CRC curtailment decisions to the network operator's actions:
\begin{align}
    &s^{\rm{el}}_{t} \leq  s_{t}, \quad \forall t \in \mathcal{T}, \label{ccc} \\
&s^{\rm{el+}}_{t} \leq  s^{+}_{t}, \quad \forall t \in \mathcal{T}, \\
& s^{\rm{el}}_{t} + s^{\rm{el+}}_{t} \leq \bar{p}^{\rm{el}}_{1} + \bar{p}^{\rm{el}}_{2}, \quad \forall t \in \mathcal{T}\label{eq: lower level CM incentivize}.
\end{align}

Constraint \eqref{eq: lower level CM incentivize} ensures that the network operator can incentivize the electrolyzer owner to opt for NFA85 or FA contracts by offering revenue from CRC activations. In  \eqref{ely_power_constraint}, \( \bar{p}^{\rm{el}}_{c} \)  is replaced by  \( \bar{p}^{\rm{no}}_{c} \) as:
\begin{equation}
     p_{t}^{\rm{e}} \leq  - s_{t} - s^{+}_{t} - r_{2,t} - r_{3,t} + \sum_{c\in{\mathcal{C}}}\bar{p}^{\rm{no}}_{c}, \quad \forall t \in \mathcal{T}, \label{LL ely_power_constraint}
\end{equation}
which ensures that the network operator cannot influence contract decisions or the minimum grid connection capacity through their curtailment actions but instead provides the available residual capacity in each hour $t$.

The following constraint is included to ensure that the electrolyzer owner has access to the residual load profile of the network, enabling them to choose contracts accordingly:
\begin{equation}
     - s^{\rm{el}}_{t} - s^{\rm{el}\mathrm{+}}_{t} - r^{\rm{el}}_{2,t} - r^{\rm{el}}_{3,t}  + \sum_{c\in{\mathcal{C}}}\bar{p}^{\rm{el}}_{c} \leq \bar{S}_{t}, \quad \forall t\in \mathcal{T}, \label{LL client grid capacity}
\end{equation}
where \( r^{\rm{el}}_{2,t} \) and \( r^{\rm{el}}_{3,t} \) represent the electrolyzer owner's curtailment decisions under the NFA85 and NFA contracts, respectively, which are restricted to the network operator's decisions:
\begin{align}
    & r^{\rm{el}}_{3,t}  \leq r_{3, t} \label{eq: ely upper nfa}, \quad \forall t \in \mathcal{T}, \\ 
    & r^{\rm{el}}_{2,t}  \leq r_{2, t} \label{eq: ely upper nfa85}, \quad \forall t \in \mathcal{T}.
\end{align}
\end{subequations}

The resulting bilevel program for the second game is therefore formulated as:
\begin{subequations} \label{bilevel NE}
\begin{align}
    & \min_{\Theta^{\rm{no,up}}} \quad \eqref{1DSO_obj} \\ 
    & \rm{s.t.} \quad \eqref{DSO_contract_constraint} - \eqref{total curtailment constraint}\ \\
    &  \bar{p}_{c}^{\rm{el}} \in  \argmin_{\Theta^{\rm{el,low}}} \{\eqref{2ely_obj_v1} \;\: \rm{s.t.} \;\:  \eqref{ely_total_contract_constraints}, \; \eqref{ely_min_load_constraint} - \eqref{1ely_upper_h2_demand}, \; \eqref{ccc} - \eqref{eq: ely upper nfa85} \}.
\end{align}
\end{subequations}

Note that all variable sets \(\Theta^{(.)}\) in both bilevel programs \eqref{bilevel EN} and \eqref{bilevel NE} belong to \(\mathbb{R}^+\).


\vspace{2mm}
\subsection{Cases}
We solve the following four cases:
\begin{enumerate}
    \item [(\textit{i})] \textbf{Game I (EL-NO)}: We solve the bilevel problem \eqref{bilevel EN}. The binary variable \( b_t \) in the lower-level problem is relaxed to a continuous variable in the range \([0, 1]\), making the lower-level problem linear and convex. We reformulate the lower-level problem using its Karush-Kuhn-Tucker optimality conditions and linearize the complementarity constraints by employing SOS1 constraints instead of Big-M formulations, ensuring optimality. To enhance tractability, we incorporate a constraint that enforces strong duality in the lower-level problem, as proposed in \cite{Kleinert2023}. 
\item [(\textit{ii})] \textbf{Game II (NO-EL)}: We solve the bilevel problem \eqref{bilevel NE}, where the lower-level problem is linear, eliminating the need for any binary relaxation. The bilevel problem is reformulated and solved in a manner similar to Game I.
\item [(\textit{iii})] \textbf{Ely HPR}: We solve the bilevel problem \eqref{bilevel EN} to its high-point relaxation (HPR). Specifically, we reformulate the bilevel problem \eqref{bilevel EN} by eliminating the lower-level objective function (but not the constraints) of the network operator. As a result, the problem becomes the optimization problem of the electrolyzer owner, additionally constrained by the feasible space of the network operator. The following constraints are also included:  
\begin{subequations}
\begin{align}
& \sum_{t \in \mathcal{T}} b_{t} \geq b^{\rm{s}} B^{\rm{NFA85}} |\mathcal{T}|, \label{eq:simplified_curtailment_85} \\
& b^{\rm{s}} \geq \frac{\bar{p}^{\rm{el}}_{2}}{M}, \label{eq:switch}
\end{align}
\end{subequations}
where $b_{t}$ and  \( b^{\rm{s}} \) are both binary variables, and \( M \)  is a sufficiently large constant. Constraints \eqref{eq:simplified_curtailment_85} and \eqref{eq:switch} prioritize the activation of NFA85 over CRC, representing the simplified curtailment decisions of the network operator. The resulting single-level optimization problem is mixed-integer linear.

\item [(\textit{iv})] \textbf{NO HPR}: We solve the bilevel problem \eqref{bilevel NE} to its HPR. Similarly, the bilevel problem \eqref{bilevel NE} is reformulated by removing the lower-level objective function of the electrolyzer owner. Consequently, the problem reduces to the optimization problem of the network operator, additionally constrained by the feasible space of the electrolyzer owner. The resulting single-level optimization problem is, once again, mixed-integer linear.
\end{enumerate}

All the above four problems are solved using Gurobi on a standard laptop, achieving a minimum optimality gap of \(0.1\%\).


\vspace{3mm}
\section{Numerical results} \label{sec: results}
Table \ref{tab:parameters} summarizes the parameter values used in the case study. Section \ref{sec: main_results} presents the contracting and grid connection capacity decisions for the different model formulations and their sensitivity to variations in \(\lambda^{\rm{CRC+}}\). A maximum hydrogen demand of 1800 kg/h is assumed, corresponding to the rated output of a 100 MW electrolyzer operating at 70\% efficiency. This capacity reflects the potential upscaling of the electrolyzer, assuming no limitations from grid congestion. The residual congestion management budget proportion, \(\theta\), is set to zero, indicating that the budget has already been fully utilized, as is often the case in regions with stalled hydrogen projects. A hydrogen selling price of 10 \texteuro/kg is assumed. Section \ref{sec: sensitivity to hydrogen price and CM budget}  analyzes the sensitivity of contracting decisions to both the hydrogen price and \(\theta\).

\begin{table}[t]
\scriptsize
\centering
\caption{Parameters}
\label{tab:parameters}
\begin{tabular}{l l l l }
\hline
\textbf{Category}       & \textbf{Parameter}       & \textbf{Value}               & \textbf{Unit} \\ \hline
\multirow{6}{*}{Electrolyzer} 
                        & Discount rate            & 0.1                          & -             \\ 
                        & Lifetime                & 15 \cite{Buttler2017}                           & years         \\ 
                        & Capital recovery factor & 0.1315                       & -             \\ 
                        & Capital cost                    & 937.5 \cite{Buttler2017}                        & k\texteuro/MW      \\ 
                        & $C^{\rm{el}}$                 & 123.26                       & k\texteuro/MW      \\ 
                        & $\alpha^{\rm{min}}$     & 0.2 \cite{Buttler2017}                        & -             \\ 
                        & $\eta^{\rm{sys}}$     & 0.7 \cite{Buttler2017}                          & -             \\ \hline
\multirow{7}{*}{Grid}   
                        & $T_{1}$                 & 87.6                         & k\texteuro/MW      \\ 
                        & $T_{2}$                 & 43.8                         & k\texteuro/MW      \\ 
                        & $T_{3}$                 & 26.28                        & k\texteuro/MW      \\ 
                        & $\bar{S}_{t}$           & time series                   & MWh             \\ 
                        & $\theta$                & \{0, 0.2\}     & -             \\ 
                        & $B^{\rm{NFA}}$          & 1                            & -             \\ 
                        & $B^{\rm{NFA85}}$        & 0.15                         & -             \\ \hline
\multirow{4}{*}{Prices} 
                        & $\lambda_{t}^{\rm{e}}$  & time series                   & \texteuro/MWh             \\ 
                        & $\lambda^{\rm{crc+}}$    & \{1-55\}        & \texteuro/MW  \\ 
                        & $\lambda^{\rm{crc}}$   & 40                           & \texteuro/MW  \\ 
                        & $\lambda^{\rm{H_{2}}}$  & \{5, 10\}       & \texteuro/kg  \\ \hline
\end{tabular}
    \vspace{-1mm}
\end{table}

\begin{figure*}[h]
    \centering
\includegraphics{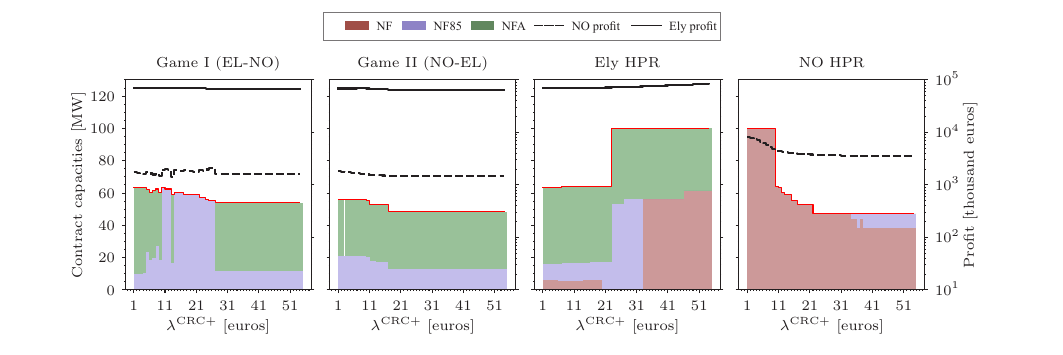}
    \caption{Optimal grid connection and CTA capacities are shown against CRC+ prices. The profits of the electrolyzer owner (Ely) and network operator (NO) are also displayed in the same figure, using a logarithmic scale.}
    \label{fig:all_formulations}
    \vspace{-1mm}
\end{figure*}


\vspace{1mm}
\subsection{Optimal Contracting Strategies} \label{sec: main_results}
Figure \ref{fig:all_formulations} presents the grid connection and CTA decisions as a function of CRC+ prices for the Game I (EL-NO), Game II (NO-EL), Ely HPR, and NO HPR cases. In the EL-NO game, the electrolyzer owner determines the CTA capacities required to secure the desired grid connection from the network operator. Recall that \( \theta = 0 \), meaning the network operator is not required to spend on activating the electrolyzer's CRC. At low CRC+ prices, the optimal electrolyzer capacity matches the network's maximum residual capacity of 63 MW. The electrolyzer owner selects contracts with the lowest tariffs up to 63 MW, first choosing 40 MW of NFA until \( B^{\rm{NFA}} \) is reached, followed by 23 MW of NFA85, which can be curtailed via CRC+ at near-zero costs. As CRC+ prices rise, the electrolyzer owner opts for larger NFA85 contracts, ensuring that the network operator's tariff income exceeds the total CRC+ costs, thus securing the desired connection capacity. As the CRC+ price increases further, the total capacity decreases to about 55 MW, at which point CRC+ activation is no longer necessary, as curtailment via NFA85 and NFA suffices.

\begin{figure}[b]
    \centering
    \includegraphics[trim={0mm 3mm 0mm 3mm},clip]{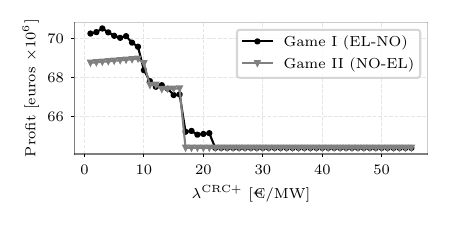}
    \caption{Profit of the electrolyzer owner in both  EL-NO and NO-EL games with the increasing CRC+ price.}
    \label{fig:profits}
\end{figure}

In the NO-EL game, the network operator takes the lead. As CRC+ prices rise, the network operator reduces the grid connection capacity to minimize costs, and the share of NFA85 contracts decreases. The electrolyzer owner does not select the FA contract. By comparing the subplots for Game I (EL-NO) and Game II (NO-EL) in Figure \ref{fig:all_formulations}, it is evident that the network operator benefits from being reactive to the electrolyzer owner's decisions at CRC+ prices above \texteuro10/MW. At these prices, the EL-NO game results in higher profits for the network operator than the NO-EL game. Figure \ref{fig:profits} zooms in on the profit of the electrolyzer owner in both games EL-NO and NO-EL. In the NO-EL game, the electrolyzer's profits increase with CRC+ prices up to \texteuro10/MW, after which they decline as the network operator reduces CRC+ activations. In contrast, in the EL-NO game, the electrolyzer owner's profit exceeds that in the NO-EL game for CRC+ prices between \texteuro1/MW and \texteuro10/MW but drops sharply as more NFA85 contracts are adopted. Higher CRC+ prices are detrimental to the electrolyzer owner, as they effectively return congestion revenue to the network operator, plus a surcharge, by choosing higher-tariff contracts, ultimately benefiting the network operator.

As expected, the HPRs of the EL-NO and NO-EL games result in larger profits for both the electrolyzer owner and the network operator, as observed  in Figure \ref{fig:all_formulations}. In the Ely HPR, it is not the network operator, but the electrolyzer owner who decides on the activation of CRC contracts. At low CRC+ prices, low-tariff contracts are preferred. As the CRC+ price increases, a trade-off emerges between tariff income and congestion costs. At CRC+ prices over \texteuro$33$/MW, the electrolyzer owner switches from NFA85 to an FA contract to benefit from additional curtailment through CRC+. This scenario reflects a current dilemma in the Netherlands, where some flexibility providers intentionally create congestion by requesting large grid connection capacities to take advantage of CRC curtailment benefits, provided a budget \( B^{\rm{CM}} \) is available for the region. However, if no budget \( B^{\rm{CM}} \) is available, this strategy cannot be employed.

In conclusion, a comparison of the NO-EL and EL-NO games with their HPRs reveals that neglecting the hierarchical structure and expanded decision space introduced by the CTA and CRC contracts leads to an overestimation of profits for both the network operator and the electrolyzer owner.

\vspace{2mm}
\subsection{Sensitivity Analysis} \label{sec: sensitivity to hydrogen price and CM budget}

\begin{figure}[t]
    \centering
    \includegraphics[trim={0mm 0mm 0mm 0mm},clip, width=0.87\linewidth]{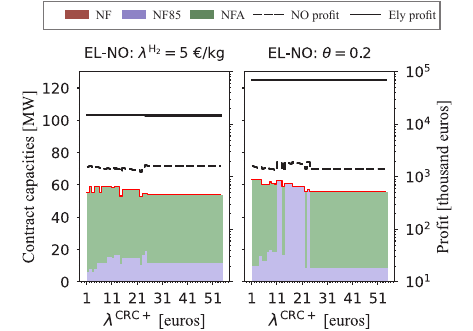}
    \caption{Game II (NO-EL): Contracting decisions at low hydrogen price (left) and 20\% congestion management budget availability (right).}
    \label{fig:sensitivtyh2theta}
\end{figure}

Figure \ref{fig:sensitivtyh2theta} shows the contracting decisions and connection capacity for the NO-EL game, assuming a hydrogen price of \texteuro5/kg, reflecting the maximum willingness to pay by potential off-takers in the GROWH project. A trade-off between grid tariffs and CRC+ prices is observed: At moderate prices, the electrolyzer owner increases the share of NFA85, but reduces it again as CRC+ prices increase further. Comparing Figures \ref{fig:all_formulations} and \ref{fig:sensitivtyh2theta}, it is evident that lowering the hydrogen selling price from \texteuro10/kg to \texteuro5/kg discourages the electrolyzer owner from choosing a full NFA85 contract, as its high tariffs render the project financially unviable. The tariff discount on NFAs helps make electrolyzer projects with uncertain transport capacity viable at low-hydrogen prices. However, Figures \ref{fig:all_formulations} and \ref{fig:sensitivtyh2theta} show that lowering the hydrogen selling price reduces grid connection capacity, as additional production fails to offset higher capital costs and tariffs.

Figure \ref{fig:sensitivtyh2theta} also shows the contracting decision with a curtailment budget availability $\theta\!=\!0.2$ and a hydrogen price of \texteuro10/kg. When the network operator is required to spend on congestion management for a larger grid connection, its influence on the electrolyzer's contract choices diminishes. This reduces the likelihood of the electrolyzer opting for the NFA85 contract, leading to lower profits for the network operator compared to the case with \( \theta = 0 \).




\vspace{2mm}
\section{conclusion} \label{sec: con}

This study investigates the interaction between an expanding electrolyzer project and a network operator in the Netherlands under new CTAs and CRCs. It formulates two bilevel programs with alternating leader-follower roles between the electrolyzer owner and the network operator. A CRC enhances electrolyzer profits at low-CRC prices (up to \texteuro10/MW), but becomes detrimental at higher prices, as the network operator is unwilling to accept low-tariff contracts while covering congestion management costs. The network operator benefits from responding to the electrolyzer owner's CTA decisions at CRC prices between \texteuro10/MW and \texteuro25/MW, where the electrolyzer is willing to pay higher tariffs to secure larger grid connection capacities, offsetting CRC costs. Non-firm contracts can make electrolyzer projects with uncertain transport capacity bankable, even a hydrogen price of \texteuro5/kg. Naive approaches that ignore the decisions of the other party overestimate the profits of both the network operator and the electrolyzer owner, underscoring the importance of a game-theoretic approach to accurately model their interaction under the new CTA and CRC contracts.


\begin{thebibliography}{11}

\bibitem{Netbeheer2025}
Netbeheer Nederland, “capaciteitkaart.” Available: \url{https://capaciteitskaart.netbeheernederland.nl/}. Accessed January 25, 2025.

\bibitem{Zoelen2024}
R. van Zoelen, N. Dooley, and D. Boer, “The role of standalone hydrogen areas in decentral hydrogen infrastructure development,” \textit{HyDelta Technical Report}, pp. 1–81, 2024. Available: \url{https://zenodo.org/records/10868624}. Accessed January 25, 2025.

\bibitem{Energiesysteem2023}
“Het energiesysteem van de toekomst: de II3050-scenario’s,” 2023. Available: \url{https://open.overheid.nl/documenten/ronl-7219ac2558977a6050ac4db764d2ddebb156df32/pdf}. Accessed January 25, 2025.

\bibitem{ACMCode2025}
“Codebesluit congestiemanagement $|$ ACM.nl.” Available: \url{https://www.acm.nl/nl/publicaties/codebesluit-congestiemanagement}. Accessed January 25, 2025.


\bibitem{Verhoeven2024}
B. van der Holst, G. Verhoeven, P. H. Nguyen, J. Morren, and K. Kok, “The Activation of congestion service contracts for budget-constrained congestion management,” \textit{Electric Power Systems Research}, vol. 235, pp. 110664, 2024.

\bibitem{Shen2018}
F. Shen, S. Huang, Q. Wu, S. Repo, Y. Xu, and J. Østergaard, “Comprehensive congestion management for distribution networks based on dynamic tariff, reconfiguration, and re-profiling product,” \textit{IEEE Transactions on Smart Grid}, vol. 10, no. 5, pp. 4795–4805, 2018.

\bibitem{Muller2023}
L. Muller and F. Cadoux, “Non-Firm grid connections for low-voltage generators: A case study in France,” in \textit{the Proceedings of the 27th International Conference on Electricity Distribution (CIRED)}, Rome, Italy, 2023.

\bibitem{Mehmood2024}
K. K. Mehmood, R. A. Arruda Moura, A. Van Der Molen, R. Tonkoski, P. Van Der Wielen, and P. H. Nguyen, “Optimization for non-firm capacity contracts for maximized grid utilization,” in \textit{the Proceedings of the 7th International Conference on Electric Power and Energy Conversion Systems}, Sharjah, United Arab Emirates, 2024.


\bibitem{ACM2022}
Autoriteit consument en markt, “Consultatie alternatieve transportrechten: use it or lose it,” 2022. Available: \url{https://www.acm.nl/nl/publicaties/consultatie-alternatieve-transportrechten-en-use-it-or-lose-it}. Accessed January 25, 2025.

\bibitem{Kleinert2023}
T. Kleinert and M. Schmidt, “Why there is no need to use a big-M in linear bilevel optimization: A computational study of two ready-to-use approaches,” \textit{Computational Management Science}, vol. 20, no. 3, pp. 1–12, 2023.

\bibitem{Buttler2017} 
A. Buttler and H. Spliethoff, “Current status of water electrolysis for
energy storage, grid balancing and sector coupling via power-to-gas and
power-to-liquids: A review,” \textit{Renewable and Sustainable Energy Reviews}, vol. 82, no. 3, pp.
2440–2454, 2018.

\end{thebibliography}
\end{document}